\documentclass[a4paper, oneside, 11pt]{article}

\usepackage{a4wide}
\usepackage{graphicx}
\usepackage{pb-diagram, lamsarrow,pb-lams}
\usepackage{amssymb}
\usepackage{parskip}
\usepackage{theorem}
\usepackage{enumerate}
\usepackage{floatflt}
\usepackage[centertags]{amsmath}
\usepackage[latin1]{inputenc}

\input xy
\xyoption{all} \CompileMatrices

\newcommand{\Hom}{\mbox{\upshape{Hom}}}

\renewcommand{\ker}{\mbox{\upshape{ker}}}
\renewcommand{\dim}{\mbox{\upshape{dim}}}
\newcommand{\Det}{\mbox{\upshape{Det}}}
\newcommand{\Qu}{\mbox{\upshape{Quot}}}
\newcommand{\res}{\mbox{\upshape{res}}}
\newcommand{\cd}{\mbox{\upshape{cd}}}

\newcommand{\Q}{\mathbb{Q}}
\newcommand{\N}{\mathbb{N}}

\newcommand{\Z}{\mathbb{Z}}

\newcommand{\QQ}{\mathcal{Q}}

\newcommand{\End}{\mbox{\upshape{End}}}

\newcommand{\tr}{\mbox{\upshape{tr}}}

\newcommand{\NN}{\mbox{\upshape{N}}}
\newcommand{\Br}{\mbox{\upshape{Br}}}
\newcommand{\nr}{\mbox{\upshape{nr}}}

\newcommand{\id}{\mbox{\upshape{id}}}

\newtheorem{St}{Step}

\newtheorem{lemma}{Lemma}
\newtheorem{theorem}{Theorem}
\newtheorem{Prop}{Proposition}
\newtheorem{corollary}{Corollary}
\newtheorem{remark}{Remark}
\newtheorem{Def}{Definition}

\begin{document}

\title{Thoughts on the reduced Whitehead group\\ of the Iwasawa algebra}

\author{Irene Lau}\date{}
\maketitle

\begin{abstract}
Let $l$ be an odd prime and $K/k$ a Galois extension of totally real number fields with Galois group $G$ such that $K/k_\infty$ and $k/\Q$ are finite. We reduce the conjectured triviality of the reduced Whitehead group $SK_1(\QQ G)$ of the algebra $\QQ G=\Qu(\Lambda G)$ with the Iwasawa algebra $\Lambda G = \Z_l[[G]]$  to the case of pro-$l$ Galois groups $G$ and finite unramified coefficient extensions.
\end{abstract}

\section{Introduction}
\label{intro}
We fix an an odd prime number $l$ and a Galois extension $K/k$ of totally real fields with Galois group $G$ such that $k/\Q$ and $K/k_\infty$ are finite.  As usual, $k_\infty$  denotes the cyclotomic $\Z_l$-extension of $k$. 
Next, the Iwasawa algebra $\Lambda G=\Z_l[[G]]=\varprojlim_{N\lhd G}\Z_l[G/N]$, where $N$ runs through the open normal subgroups of $G$, denotes the completed group ring of $G$ over $\Z_l$ and $\QQ G=\Qu(\Z_l[[G]])$ is its total ring of fractions with respect to all central non-zero divisors. Let $K_0T(\Lambda G)$ be the Grothendieck group of the category of finitely generated torsion $\Lambda G$-modules of finite projective dimension. Then, the localization sequence of $K$-theory
$$\rightarrow K_1(\Lambda G)\rightarrow K_1(\QQ G)\stackrel{\partial}{\rightarrow}K_0T(\Lambda G)\rightarrow$$
is exact.
$\QQ G$ finds its way into non-commutative Iwasawa theory via this localization sequence and a determinant map
$$\Det: K_1(\QQ G)\rightarrow \Hom(R_lG, (\Q_l^c\otimes_{\Q_l}\QQ \Gamma_k)^\times),$$
where $\Q_l^c$ is a fixed algebraic closure of $\Q_l$, the $\Z_l$-span of the irreducible $\Q_l^c$-characters of $G$ with open kernel is named $R_lG$ and $\Gamma_k=G(k_\infty/k)$. This determinant is the translation of the reduced norm $\nr: K_1(\QQ G)\rightarrow Z(\QQ G)^\times$ to Hom groups, where $Z(\QQ G)$ is the centre of $\QQ G$.
We refer to \cite{Iw2} for a precise definition of $\Det$.

As in the classical case of Iwasawa,  Ritter and Weiss link a $K$-theoretic substitute $\mho\in K_0T(\Lambda G)$ of the Iwasawa module $X$ to the Iwasawa $L$-function which is derived from the $S$-truncated Artin $L$-function for a finite set $S$ of places of $k$ containing all archimedian ones and those which ramify in $K$ (see e.g. \cite{Iw2}). This Iwasawa $L$-function lies in the upper $\Hom$-group. 

With this, the main conjecture of equivariant Iwasawa theory says
$$There\ exists\ a\ unique\ element\ \Theta\in K_1(\QQ G)\ s.t.\ \Det(\Theta)=L.\ Moreover\ \partial(\Theta)=\mho.$$
The uniqueness of $\Theta$ would follow from the conjecture by Suslin 
that the reduced Whitehead group $SK_1(A)$ is trivial for central simple algebras $A$ over fields with cohomological dimension $\leq 3$ (see \cite{Sus} for this conjecture of Suslin). In \cite[Thm 5.1]{QG}, it is shown that this conjecture can be applied to our algebra $\QQ G$.

Recently (compare \cite{Iw10}), Ritter and Weiss gave a complete proof of this main conjecture up to its uniqueness statement whenever Iwasawa's $\mu$-invariant vanishes. In \cite{Kak2}, Kakde also gave a proof. In fact, he does not restrict to 1-dimensional $l$-adic Lie groups as Ritter and Weiss do but gives a proof for higher dimensional admissible $l$-adic Lie groups. Yet, he does not consider the full ring of fractions $\QQ G$ but the localization $(\Lambda G)_S$ by the canonical Ore set $S$ and proves uniqueness up to the quotient of $K_1((\Lambda G)_S)$ by the image of $SK_1(\Lambda G))$. Thus, the question whether $\Theta$ is unique in $K_1(\QQ G)$ is still open.

In this paper, we reduce the Suslin conjecture for our Iwasawa algebra $\QQ G$ for profinite Galois groups $G$ to the conjecture for $N\otimes_{\Q_l}\QQ U$ for pro-$l$ groups $U$ and finite unramified extensions $N$ of $\Q_l$. Therefore, the proof of the uniqueness statement of the main conjecture is completely reduced to pro-$l$ groups provided that the studied objects are unaffected by passing to finite unramified extensions of $\Q_l$.

This paper contains some of the results of my PhD thesis. I would like to thank my supervisor Jürgen Ritter for his aid, encouragement and patience during my work on this paper.

\section{Recollections}\label{Recollections}
First, we recall some facts on the structure of $\QQ G$ and formulate the Suslin conjecture.

We keep the notation of the introduction, in particular we fix an odd prime $l$ and a Galois extension $K/k$ of totally real fields with Galois group $G$ such that $k/\Q$ and $K/k_\infty$ are finite.

First, $G$ splits (see \cite[p. 551]{Iw2}): 
$G=H\rtimes \Gamma$ with $H=G(K/k_\infty)$ and $\Gamma=\langle\gamma\rangle\cong G(k_\infty/k)\cong \Z_l$.
 Thus, for a central subgroup $\Gamma^{l^m}=:\Gamma_0$ we get  $$\QQ G=\bigoplus_{i=0}^{l^m-1}(\QQ\Gamma_0)[H]\gamma^{i}.$$
This algebra is a finite dimensional $\QQ\Gamma_0$-algebra; in fact, it is a semisimple algebra, since 
 the Jacobson radical is trivial by \cite[p. 553]{Iw2}.
  Now, let $\chi\in R_lG$ be an irreducible $\Q_l^c$-character of $G$ with open kernel. Note that it is sufficient to regard the finite set of irreducible characters of $G/\Gamma_0$ because, by inflation and twist with irreducible characters $\rho$ which fulfil $\res_G^H\rho=1$, every irreducible $\chi\in R_lG$ can be obtained from this set. These characters $\rho$ will be called of type W. Because $G$ is an $l$-group with $l\neq 2$, this implies that $\chi$ has a representation over $\Q_l(\chi)$ by \cite{Ro}. 

 Furthermore, with $\eta$  an absolutely irreducible constituent of $\res_G^H(\chi)$, we define 
$$St(\eta):=\{g\in G:\eta^g=\eta\},\ w_\chi:=[G:St(\eta)]$$
 and 
 $$e(\eta):=\frac{\eta(1)}{|H|}\sum_{h\in H}\eta(h^{-1})h.$$
Ritter and Weiss showed in \cite{Iw2} that 
$$e_\chi:=\sum_{\eta|\mbox{\scriptsize{res}}_G^H\chi}e(\eta)$$
 is a central primitive idempotent in $\QQ^cG$, that every  central primitive idempotent is of the form $e_\chi$ and that two  central primitive idempotents $e_{\chi_1}$ and $e_{\chi_2}$ coincide if and only if $\chi_1=\chi_2\otimes \rho$ for a character $\rho$ of type $W$.

In the special case of a pro-$l$ group $G$, the structure of $\QQ G$ is completely known by \cite{QG}: 
\begin{lemma}\label{lemQG}
Let now $G$ be a pro-$l$ group and let $W'$ be the simple component of $\QQ G$ corresponding to the irreducible character $\chi\in R_l G$. We moreover choose an absolutely irreducible constituent $\eta$ of $\res_G^H(\chi)$. Then, 
\begin{enumerate}
\item $$W'\cong \bigoplus_{i=0}^{l^m-1}\left(\bigoplus_{j=0}^{v_\chi-1}(\Q_l(\eta)\otimes_{\Q_l}\QQ\Gamma_0)_{\eta(1)\times\eta(1)}\right)\gamma^{i}
$$
 for $v_\chi:=\mbox{min}\{0\leq j\leq w_\chi-1 : \eta^{\gamma^j}=\eta^\sigma \mbox{ for some } \sigma\in G(\Q_l(\eta)/\Q_l)\}$,
\item $W'$ has centre 
$$Z(W')\cong L\otimes_{\Q_l}\QQ\Gamma^{w_\chi}$$
with $L=\Q_l(\eta)^{G_0}$ and $G_0=\{\sigma\in G(\Q_l(\eta)/\Q_l):\eta^\sigma=\eta^{\gamma^j} \mbox{ for a } 0\leq j\leq w_\chi-1\}$.\\
Moreover, $G_0=:\langle \sigma_{v_\chi}\rangle$ is a cyclic group of order $\frac{w_\chi}{v_\chi}$.
\item $Z(W')$ has cohomological dimension $\cd(Z(A))=3$,
\item $\dim_{Z(W')}W'=\chi(1)^2$,
\item $W'$ is split by $\Q_l(\eta)\otimes_{\Q_l}\QQ\Gamma^{w_\chi}$, 
\item $W'$ has Schur index $s_D=w_\chi/v_\chi$ and 
\item $W'\cong D_{n\times n}$ with $n=\chi(1)/s_D$ and the skew field $D$ is cyclic:
$$D\cong  \bigoplus_{i=0}^{w_\chi/v_\chi-1} (\Q_l(\eta)\otimes_{\Q_l}\QQ\Gamma^{w_\chi}) \gamma^{v_\chi i}=: (\Q_l(\eta)\otimes_{\Q_l}\QQ\Gamma^{w_\chi}/L\otimes_{\Q_l}\QQ\Gamma^{w_\chi},\sigma_{v_\chi},\gamma^{w_\chi}).$$
\end{enumerate}
\end{lemma}
\textbf{Proof:}
Statements (i) and (ii) can be found in \cite[Prop 1]{QG}, (iii) is \cite[Thm 2]{QG} and \cite[Thm 1]{QG} contains (iv) to (vii).\hfill$\Box$

Because $H$ is a finite $l$-group, $\Q_l(\eta)$ is generated by a primitive $l$-power root of unity. Therefore, $L=\Q_l(\eta)^{G_0}\subseteq \Q_l(\eta)$ also is, i.e. we can fix a primitive $l$-power root of unity $\xi$ s.t.
$$L=\Q_l(\xi).$$

We now focus on the Suslin conjecture.
\begin{Def}
\begin{enumerate}
\item For a field $F$  and a central simple $F$-algebra $A$ of finite degree $[A:F]$, let $\nr_{A/K}$ denote the reduced norm from $A$ to $K$. The group
$$SK_1(A):=\ker(\nr_{A/F})/[A^\times,A^\times]$$
is called the reduced Whitehead group of $A$.
\item For a semisimple algebra $A=\bigoplus_i A_i$ of finite degree with simple components $A_i$, we set
$$SK_1(A):=\bigoplus_i SK_1(A_i)$$
for the reduced Whitehead group of $A$.
\end{enumerate}
\end{Def}

The reduced norm $\nr_{A/F}$ on $A$ induces a homomorphism on $K_1(A)$, which we will call reduced norm, too. We state the following well-known results without proof.
\begin{lemma}\label{1}
\begin{enumerate}
\item Let $A$ be a central simple $F$-algebra of finite degree. Then
$$SK_1(A)=\ker(\nr_{A/F}:K_1(A)\rightarrow K_1(F)).$$
\item Let $A\cong D_{n\times n}$ be the full matrix ring of finite degree over a skew field $D$. Then
$$SK_1(A)=SK_1(D).$$
\item For a field $F$, we have
$$SK_1(F)=1.$$
\end{enumerate}
\end{lemma}
For further details, see e.g.~\cite[Part III]{Drax}. 

\begin{remark}
The determonant map $\Det$ in the main conjecture of equivariant Iwasawa theory is the translation of the reduced norm to the language of Hom-groups. For a detailed definition of this $Det$, we refer to \cite[p. 558]{Iw2}.
\end{remark}

We are now ready to state the

\textbf{Conjecture}
\textit{Let $F$ be a field with cohomological dimension $\cd(F)\leq 3$ and $A$ a central simple $F$-algebra of finite degree $[A:F]$. Then
$$SK_1(A)=1.$$
}

In the following, we will call this  Suslin's conjecture, although this is not literally  Suslin's formulation. But in the case of a field of cohomological dimension less than or equal to $3$, this is exactly the statement of his conjecture. For details, we refer to \cite{Sus}.

The centres of the Wedderburn components of $\QQ G$, i.e. the simple components $W'$, are of cohomological dimension $3$ for pro-$l$ groups $G$ by Lemma \ref{lemQG}. As we will see in this paper, this is the crucial case for the triviality of $SK_1(\QQ G)$.

Next, we list the cases for $\QQ G$ which are known to have trivial reduced Whitehead group:
\begin{lemma}
Let $G$ be as above. Then, $SK_1(\QQ G)=1$ in the following cases:
\begin{enumerate}
\item $G=H\times \Gamma$ is a direct product with $H$ an $l$-group or of order prime to $l$\footnote{We will see in subsection \ref{l-q-section} that the restrictions on $H$ are not necessary.}.
\item $G$ is a pro-$l$ group $G$ with abelian subgroup of index $l$.
\item $G=H\rtimes \Gamma$, where $H$ is a finite group of  order prime to $l$.
\end{enumerate} 
\end{lemma}
\textbf{Proof:}
For $l$-groups $H$, Roquette has shown in \cite{Ro} that $\Q_l[H]$ is the direct sum of some matrix rings over fields.
Therefore, $\QQ G=(\QQ \Gamma)[H]=\QQ \Gamma\otimes_{\Q_l}\Q_l[H]$ also is a direct sum of matrix rings over fields and thus $SK_1(\QQ G)=1$ by Lemma \ref{1}. This shows the first case of (i). The latter statement of (i) is a special case of (iii).

(ii) is shown in \cite[p.~118]{wall}.

(iii) can be found in \cite[Example 2, p.~169]{Iw4}.\hfill$\Box$

\section{Reduction of the reduced Whitehead group $SK_1(\QQ G)$ to the pro-$l$ case}\label{reduction2}
 
As main ingredient of our reduction, we cite the following lemma (see \cite[Cor, p.~167]{Iw4}). For this, recall that, for a prime number $q\neq l$, $G$ is  a $\Q_l$-$q$-elementary  group if $G=H\times \Gamma$ with $\Gamma$ a central open subgroup of $G$ isomorphic to $G(k_\infty/k)$ 
 and $H$ a finite $\Q_l$-$q$-elementary group; i.e.~$H=\langle s \rangle \rtimes H_q$ is the semidirect product of a cyclic group $\langle s \rangle$ of order prime to $q$ and a $q$-group $H_q$ whose action on $\langle s \rangle$ induces a homomorphism $H_q\rightarrow G(\Q_l(\zeta)/\Q_l)$. Here, $\zeta$ is a primitive root of unity of order $|\langle s\rangle |$. For $q=l$, the group $G$ is called $\Q_l$-$l$-elementary if $G=\langle s\rangle\rtimes U$ is the semidirect product of a finite cyclic group $\langle s \rangle$ of order prime to $l$ and an open pro-$l$ subgroup $U$ whose action on $\langle s \rangle$ induces a homomorphism $U\rightarrow G(\Q_l(\zeta)/\Q_l)$ with again $\zeta$ a primitive root of unity of order $|\langle s \rangle|$.
\begin{lemma}[Ritter, Weiss]
Let $K/k$ be a Galois extension of totally real fields with Galois group $G$ such that $K/k_\infty$ and $k/\Q$ are finite. Then, $SK_1(\QQ G)=1$ if $SK_1(\QQ G')=1$ for all open $\Q_l$-$q$-elementary subgroups $G'$ of $G$ and all prime numbers $q$ ($q$ might be equal to $l$). 
\end{lemma}
Thus, we have to compute $SK_1(\QQ G)$ for $\Q_l$-$q$-elementary groups $G$ with $q$ running through the set of all prime numbers.

\subsection{$\Q_l$-$l$-elementary groups $G$}
We begin with the case $q= l$, i.e.~$G=\langle s \rangle \rtimes U$ with a finite cyclic group $\langle s \rangle$ of order prime to $l$ and $U$ an open pro-$l$ subgroup. 

We fix a finite set $\{\beta_i\}$ of representatives of the $G(\Q_l^c/\Q_l)$-orbits of the irreducible $\Q_l^c$-characters of $\langle s \rangle$. Let also $\zeta_i$ denote a fixed primitive $l$-prime root of unity with $\beta_i(s)=\zeta_i$. 

Let $U_i:=\{u\in U: \beta_i^{u}=\beta_i\}$ denote the stabilizer group of $\beta_i$.  Clearly, $U_i\lhd U$ and $A_i:=U/U_i\leq G(\Q_l(\beta_i)/\Q_l)=G(\Q_l(\zeta_i)/\Q_l)$. Thus, $A_i$ is cyclic because $\Q_l(\beta_i)=\Q_l(\zeta_i)$ is unramified over $\Q_l$. We fix a representative $x_i\in U$ with $\langle \overline{x_i}\rangle =U/U_i=A_i$.  Then, $\overline{x_i}$ maps to some $\tau_i$ under the injection $U/U_i\rightarrowtail G(\Q_l(\beta_i)/\Q_l)=G(\Q_l(\zeta_i)/\Q_l)$  and therefore the order of $\overline{x_i}$ clearly is a power of $l$, say $l^n:=|U/U_i|$.
 Although $n$  depends on $i$, we omit this in the notation. Moreover, for the sake of brevity, we set $x:=x_i$ and $\tau:=\tau_i$, but still keep in mind the underlying $\beta_i$.
Finally, we set $G_i:=\langle s \rangle \rtimes U_i$.

We next read the structure of $\QQ G$ in these terms. For this, recall that the  
$$e_i:=\frac{1}{|\langle s \rangle|} \sum_{\nu\mbox{\footnotesize{\upshape{mod}}} |\langle s\rangle|} \tr_{\Q_l(\zeta_i)/\Q_l}(\zeta_i(s^{-\nu}))s^\nu\in \Z_l\langle s \rangle$$
 are the primitive central idempotents of the group algebra $\Q_l\langle s\rangle$, and furthermore they are central idempotents of $\QQ G$. Because the $e_i$ are orthogonal in $\Q_l\langle s\rangle$, we have $e_ie_j=0$ for $i\neq j$ in $\QQ G$, too. Therefore, we conclude $\bigoplus_{i}e_i\QQ G\subseteq \QQ G$. For
$$\QQ G=\bigoplus_i e_i\QQ G,$$
it remains to show the other inclusion $\QQ G\subseteq \bigoplus_{i}e_i\QQ G$. We use that $\sum_i e_i=1$ is true in $\Q_l\langle s \rangle$ and therefore it is true in $\QQ G$, too. Thus, $\QQ G=1\cdot\QQ G\subseteq \bigoplus_{i}e_i\QQ G$. We are now ready to state

\begin{lemma}\label{lem4.1}
With the above notations, we have:
\begin{enumerate}
\item $e_i \QQ G_i\cong \Q_l(\zeta_i)\otimes_{\Q_l}\QQ U_i$,
\item  $e_i \QQ G\cong \bigoplus_{j=0}^{l^n-1}(\Q_l(\zeta_i)\otimes_{\Q_l}\QQ U_i)x^j$,\\
 where $x$ acts on $U_i$ by conjugation and on $\Q_l(\zeta_i)$ via $\tau$.
\end{enumerate}
\end{lemma}
\textbf{Proof:}
(i) is stated in  \cite[p.~160]{Iw4} and (ii) follows immediately by (i) and the definition of $U_i$. \hfill$\Box$

To point out the importance of the operation of $x$, we will also use the notation
$$(\Q_l(\zeta_i)\otimes_{\Q_l}\QQ U_i)\star \langle x \rangle:=\bigoplus_{j=0}^{l^n-1}(\Q_l(\zeta_i)\otimes_{\Q_l}\QQ U_i)x^j.$$
 
\begin{Prop}
With the above notations, the following are equivalent:
\begin{enumerate}
\item $SK_1(\QQ G)=1$.
\item $SK_1(e_i \QQ G)=SK_1((\Q_l(\zeta_i)\otimes_{\Q_l}\QQ U_i)\star \langle x\rangle)=1$ for all characters $\beta_i$ of $\langle s \rangle$.
\end{enumerate}
\end{Prop}
\textbf{Proof:}
This follows immediately by $\QQ G=\bigoplus_i e_i\QQ G$ and Lemma \ref{lem4.1}.
\hfill$\Box$

As the structure of $\QQ U_i$ is well known from Section \ref{Recollections}, we  now examine $(\Q_l(\zeta_i)\otimes_{\Q_l}\QQ U_i)\star \langle x\rangle$.\\
 Because $(\Q_l(\zeta_i)\otimes_{\Q_l}\QQ U_i)\star \langle x\rangle$ is isomorphic to $e_i\QQ G$, this algebra is semisimple.

Let $W'$ be the Wedderburn component, i.e. the simple component, of $\QQ U_i$ corresponding to $\chi\in R_lU_i$ and set 
$$W=\Q_l(\zeta_i)\otimes_{\Q_l} W'=(\Q_l(\zeta_i)\otimes_{\Q_l} Z(W'))\otimes_{Z(W')} W'\subseteq \Q_l(\zeta_i)\otimes_{\Q_l}\QQ U_i.$$
As $\Q_l(\zeta_i)$ and $F':=Z(W')=L\otimes_{\Q_l} \QQ \Gamma^{w_\chi}$ are linearly disjoint over $\Q_l$, the tensor product $\Q_l(\zeta_i)\otimes_{\Q_l} F'$ is a field and thus $W$ is still a simple algebra and therefore a Wedderburn component of $\Q_l(\zeta_i)\otimes_{\Q_l}\QQ U_i$ with centre $F:=Z(W)=\Q_l(\zeta_i)\otimes F'$.

 Then, $x$ acts on $W$ as it acts on $\Q_l(\zeta_i)\otimes_{\Q_l}\QQ U_i$. This action fixes the algebra $\Q_l(\zeta_i)\otimes_{\Q_l}\QQ U_i$ as a whole, but might not fix $W$. If $W^x\neq W$, then $W^x$ is another Wedderburn component of $\Q_l(\zeta_i)\otimes_{\Q_l}\QQ U_i$ by the following: $W^x$ is a two-sided ideal of $\Q_l(\zeta_i)\otimes_{\Q_l}\QQ U_i$ because $W$ is a two-sided ideal of $\Q_l(\zeta_i)\otimes_{\Q_l}\QQ U_i$. Furthermore, it has centre $F^x$ with $F=Z(W)$. As seen above, $F=\Q_l(\zeta_i)\otimes_{\Q_l} L\otimes_{\Q_l} \QQ \Gamma^{w_\chi}$ is a field and therefore $F^x$ is a field, too. But as a semisimple algebra with a field as centre, $W^x$ is already a simple algebra. Thus, $x$ permutes the Wedderburn components of $\Q_l(\zeta_i)\otimes_{\Q_l}\QQ U_i$ and $W^x\cdot W=0$ if $W^x\neq W$ because of the orthogonality of Wedderburn components.

Note that the minimal $j$, such that $W^{x^j}=W$, is an $l$-power because this is the length of the orbit of $W$ in the set of Wedderburn components of $\Q_l(\zeta_i)\otimes_{\Q_l}\QQ U_i$ under the action of $\langle \overline{x}\rangle$.

\begin{Prop}\label{st}
Let $W$ be a simple component of $\Q_l(\zeta_i)\otimes_{\Q_l}\QQ U_i$ with centre $F$. Set $0\leq d\leq n$ to be minimal such that $W^{x^{l^d}}=W$. Then,
$$\tilde{W}:=\bigoplus_{j=0}^{l^d-1}(W^{x^j}\oplus W^{x^j}x\oplus...\oplus W^{x^j}x^{l^n-1})$$
is a simple component of $(\Q_l(\zeta_i)\otimes_{\Q_l}\QQ U_i)\star \langle x\rangle$ with centre $Z(\tilde{W})=F^{\langle x^{l^d}\rangle}=:E$.\\
Furthermore, $\tilde{W}$ is the full matrix ring 
$$\tilde{W}=V_{l^d\times l^d} \ \ \mbox{with} \ \ V:=W\oplus Wx^{l^d}\oplus...\oplus Wx^{l^d(l^{n-d}-1)}.$$
\end{Prop}
\textbf{Proof:}
We  set $y:=x^{l^d}$ and $m:=n-d$, i.e.~$y^{l^m}=x^{l^n}\in U_i$. 

First, $\tilde{W}$ is a two-sided ideal of $(\Q_l(\zeta_i)\otimes_{\Q_l}\QQ U_i)\star \langle x\rangle$; for this, we only have to check that it is closed under multiplication with $x$, which is obvious. Thus it is the direct sum of some Wedderburn components. 

Next, we show that the centre of $\tilde{W}$ is a field, which automatically implies that $\tilde{W}$ is a simple algebra. We start with the computation of the centre $Z(V)$ of $V$. Here, we do not consider the trivial case  $d=n$, i.e.  $V=W$ and therefore $Z(V)=Z(W)=F$ is a field.

We  assume $0\leq d<n$ and  take an element $z=w_0+w_1y+...+w_{l^m-1}y^{l^m-1}\in Z(V)$. For any $w\in W\subseteq V$, we see
\begin{align*}
zw&=w_0w+w_1w^{y^{-1}}y+...+w_{l^m-1}w^{y^{-(l^m-1)}}y^{l^m-1},\\
wz&=ww_0+ww_1y+...+ww_{l^m-1}y^{l^m-1}.
\end{align*}
Because $zw=wz$, we conclude that $w_0\in Z(W)=F$ and 
\begin{align}\label{z1}
w_1w^{y^{-1}}=ww_1\ ,\ ...\ ,\ w_{l^m-1}w^{y^{-(l^m-1)}}=ww_{l^m-1}.
\end{align} 

Assume for the moment that $w\in Z(W)=F$. 
 Then, (\ref{z1}) implies 
 $$w_1w^{y^{-1}}=w_1w\ ,\ ...\ ,\ w_{l^m-1}w^{y^{-(l^m-1)}}=w_{l^m-1}w.$$ But as $F=\Q_l(\zeta_i)\otimes_{\Q_l}L\otimes_{\Q_l} \QQ \Gamma^{w_\chi}$, we can specialize to $w=\zeta_i$. By definition, $y$ does not act trivially on $\zeta_i$ (otherwise $y\in U_i$) and thus $w_1=...=w_{l^m-1}=0$.

Moreover, $z$ fulfils $yz=zy$. As we have already seen that $z=w_0\in F$, this implies that $z\in F^{\langle y\rangle}$. Thus $Z(V)\subseteq F^{\langle y\rangle}$. Because the other inclusion $Z(V)\supseteq F^{\langle y\rangle}$ is trivially true, we finally conclude
$$Z(V)=F^{\langle y\rangle}.$$

Now, we are ready to show that $Z(\tilde{W})=Z(V)=F^{\langle y \rangle}$. For the rest of the proof, we will again allow the trivial case, i.e.~$0\leq d\leq n$. We use the relation
\begin{align*}
\tilde{W}&=\bigoplus_{j=0}^{l^d-1}(W^{x^j}\oplus W^{x^j}x\oplus...\oplus
W^{x^j}x^{l^n-1})\\
&=\bigoplus_{j=0}^{l^d-1}(V^{x^j}\oplus V^{x^j}x\oplus...\oplus V^{x^j}x^{l^d-1}).
\end{align*}
Let $0\leq j\leq l^d-1$. 
Because $W^{x^j}\neq W$, we have seen $W^{x^j}\cdot W=0$ and therefore $V^{x^j}\cdot V=0$.

We choose $z=\sum_{i,j=0}^{l^d-1} v_{ij}^{x^j}x^{i}\in Z(\tilde{W})$, and $v, v'\in V$. Then
\begin{align*}
zv&=\sum_{i,j}v_{ij}^{x^j}v^{x^{-i}}x^i=v_{00}v+\sum_{i>0}(v_{i,l^d-i}v^{x^{-l^d}})^{x^{l^d-i}}x^i\in V\oplus\bigoplus_{i=1}^{l^d-1}V^{x^{l^d}-i}x^i,\\
vz&=\sum_{i,j}vv_{ij}^{x^j}x^i=\sum_ivv_{i0}x^i=vv_{00}+\sum_{i>0}vv_{i0}x^i\in V\oplus\bigoplus_{i=1}^{l^d-1}Vx^i.
\end{align*}
Thus, $v_{00}\in Z(V)$ and the orthogonality of the $V^{x^j}$ implies $v_{i,l^d-i}=0=v_{i0}$ for all $i>0$. Next,
\begin{align*}
zv^x&=\sum_{i,j}v_{ij}^{x^j}v^{x^{-i+1}}x^i=(v_{01}v)^x+v_{10}vx+\sum_{i>1}(v_{i,l^d-i+1}v^{x^{-l^d}})^{x^{l^d-i+1}}x^i,\\
v^xz&=\sum_{i,j}v^xv_{ij}^{x^j}x^i=\sum_i(vv_{i1})^xx^i=(vv_{01})^x+(vv_{11})^x x+\sum_{i>1}(vv_{i1})^xx^i.
\end{align*}
Thus, $v_{01}\in Z(V)$, $v_{10}=0=v_{11}$ and $v_{i,l^d-i+1}=0=v_{i1}$ for all $i>1$. Analogous computations for $zv^{x^\nu}=v^{x^\nu}z$ finally lead to $$z=v_0+v_1^x+...+v_{l^d-1}^{x^{l^d-1}}$$
 with $v_i\in Z(V)$. We apply this together with the orthogonality and compute
\begin{align*}
z(v+v'x^i)&=v_0v+v_0v'x,\\
(v+v'x^i)z&=vv_0+v'v_ix=v_0v+v_iv'x
\end{align*}
with $1\leq i\leq l^d-1$. Thus, $v_i=v_0$ for all $1\leq i\leq l^d-1$. 

Therefore, we have achieved
$z\in\{\sum_{j=0}^{l^d-1}v^{x^j}:v\in Z(V)\}\cong Z(V)$,
i.e.~$Z(\tilde{W})\subseteq Z(V)$.
 For the other inclusion, it remains to show that elements  of $\{\sum_{j=0}^{l^d-1}v^{x^j}:v\in Z(V)\}$ are already central in $\tilde{W}$. For this, we only have to check that $\sum_{j=0}^{l^d-1}v^{x^j}$ commutes with $x$ for every $v\in Z(V)=F^{\langle x^{l^d}\rangle}$:
$$\left(\sum_{j=0}^{l^d-1}v^{x^j}\right)^x=\sum_{j=0}^{l^d-1}v^{x^{j+1}}=\sum_{j=1}^{l^d-1}v^{x^j}+v^{x^{l^d}}=\sum_{j=1}^{l^d-1}v^{x^j}+v=\sum_{j=0}^{l^d-1}v^{x^j}.$$
Hence, $Z(\tilde{W})=Z(V)$ is true. This moreover shows that $\tilde{W}$ is a Wedderburn component of $(\Q_l(\zeta_i)\otimes_{\Q_l}\QQ U_i)\star\langle x\rangle$.

It remains to show that $\tilde{W}$ is the claimed matrix ring. First, if $\tilde{W}$ is a matrix ring over $V$, then the dimension is clear because
\begin{align*} \dim_{Z(\tilde{W})}\tilde{W}&=\dim_{Z(V)}\tilde{W}=\dim_{Z(V)}\bigoplus_{j=0}^{l^d-1}(V^{x^j}\oplus V^{x^j}x\oplus...\oplus V^{x^j}x^{l^d-1})\\
&=l^{2d}\dim_{Z(V)}V
\end{align*}
and therefore $\dim_V \tilde{W}=\dim_{Z(V)}\tilde{W}/\dim_{Z(V)}V=l^{2d}$.

Both $V$ and $\tilde{W}$ are central simple algebras over $F^{\langle y\rangle}$. We show that $V\sim \tilde{W}$ in $\Br(F^{\langle y \rangle})$, i.e.~that $V$ and $\tilde{W}$ are full matrix rings over the same skew field $D$ of centre $F^{\langle y\rangle}$.

For the  computation of the skew field $D$ in $\tilde{W}$, we recall the fact that there exists a primitive idempotent $\varepsilon$ of $\tilde{W}$ such that $D\cong \varepsilon \tilde{W}\varepsilon$ and $\tilde{W}\cong B^n$   for a minimal right ideal $B=\varepsilon \tilde{W}$ of $\tilde{W}$.
Analogously, there exists a primitive idempotent $\varepsilon_V\in V$ with $\varepsilon_VV\varepsilon_V\cong D_V$ a skew field and $S=\varepsilon_V V$ a minimal right ideal of  $V$. Then, we get $\varepsilon_V \tilde{W} \varepsilon_V=\varepsilon_V V\varepsilon_V$ because for $\sum_{i,j=0}^{l^d-1}v_{ij}^{x^j}x^i\in \tilde{W}$, we achieve
\begin{align*}
\varepsilon_V\cdot \left(\sum_{i,j=0}^{l^d-1}v_{ij}^{x^j}x^i\right)\cdot \varepsilon_V
&=\varepsilon_V\cdot\left( \sum_{i,j=0}^{l^d-1}v_{ij}^{x^j}\varepsilon_V^{x^{-i}}x^i\right)\\
&\stackrel{1}{=}\varepsilon_V v_{00}\varepsilon_V+\varepsilon_V\cdot\left(\sum_{i=1}^{l^d-1}(v_{i,l^d-i}^{x^{l^d}})^{x^{-i}} \varepsilon_V^{x^{-i}}x^i\right) \\
&=\varepsilon_V v_{00} \varepsilon_V+\left(\sum_{i=1}^{l^d-1}\varepsilon_V(v_{i,l^d-i}^{x^{l^d}})^{x^{-i}}x^i\right)\cdot\varepsilon_V\\
&\stackrel{2}{=}\varepsilon_Vv_{00}\varepsilon_V.
\end{align*}
For $\stackrel{1}{=}$ and $\stackrel{2}{=}$, we have again used $V^{x^j}\cdot V=0$ for $1\leq j\leq l^d-1$.

Next, as $\varepsilon_V \tilde{W}$ is a right ideal of $\tilde{W}$, there exists a $0<r\in \N$ with $B^r\cong\varepsilon_V\cdot\tilde{W}$ for the minimal right ideal $B$. Because $\End_{\tilde{W}}(B)\cong D$ is the skew field lying in $\tilde{W}$, we get
\begin{align*}
D_V\cong\varepsilon_V V\varepsilon_V &=\varepsilon_V \tilde{W} \varepsilon_V\cong \End_{\tilde{W}}(\varepsilon_V\tilde{W})\\
&\cong\End_{\tilde{W}}(B^r)\cong\End_{\tilde{W}}(B)_{r\times r}\cong D_{r\times r}.
\end{align*}
This forces $r=1$ (because, for example, $D_V$ does not have zero divisors whereas $D_{r\times r}$ has for $r>1$).
Thus, the underlying skew fields of $V$ and $\tilde{W}$ are equal, i.e.~$V\sim\tilde{W}$ in $\Br(E)$. By $V\subseteq \tilde{W}$, this implies the claim and concludes the proof.\hfill$\Box$

\begin{remark}
The $\tilde{W}$ in Proposition \ref{st} exhaust all simple components of $(\Q_l(\zeta_i)\otimes_{\Q_l}\QQ U_i)\star\langle x\rangle$ because $$\Q_l(\zeta_i)\otimes_{\Q_l}\QQ U_i=\bigoplus_W W$$
 and therefore
$$(\Q_l(\zeta_i)\otimes_{\Q_l}\QQ U_i)\star\langle x\rangle=\left(\bigoplus_W W\right)\star \langle x\rangle=\bigoplus_{\tilde{W}}\tilde{W}.$$
\end{remark}
\begin{corollary}\label{Cor1}
With the notations of Proposition \ref{st}, we have
$$SK_1(\tilde{W})=SK_1(V).$$
\end{corollary}
\textbf{Proof:}
This is obvious, as the reduced Whitehead group of a simple algebra only depends on the underlying skew field but not on the matrix degree.\hfill$\Box$

Since the Wedderburn components $\tilde{W}$ are classified now, we can start our study of $SK_1(\tilde{W})$.

\begin{theorem}\label{sk1}
Let $G=\langle s\rangle\rtimes U$ be a $\Q_l$-$l$-elementary group with a finite cyclic group $\langle s \rangle$ of order prime to $l$ and $U$ an open pro-$l$ subgroup. Assume that $SK_1(\Q_l(\zeta_i)\otimes_{\Q_l}\QQ U_i)=1$ for all $i$. Then 
 $$SK_1(\QQ G)=1.$$
\end{theorem}

 The proof of this theorem depends on the number $l^d=\min\{1\leq j\leq l^n:\ W^{x^{j}}=W\}$.

\subsubsection{The case $d=n$}

First, let $d=n$, i.e.~$W^{x^{j}}\neq W$ for all $1\leq j\leq l^n-1$. Thus,
$$\tilde{W}=\bigoplus_{j=0}^{l^n-1} (W^{x^j}\oplus...\oplus W^{x^j}x^{l^n-1}),\ \ V=W.$$
Then, Proposition \ref{st} implies that $\tilde{W}=V_{l^n\times l^n}=W_{l^n\times l^n}$. Furthermore, both $W$ and $\tilde{W}$ have centre $Z(W)=Z(\tilde{W})=F$. (Observe that $F$ commutes with $x^{l^n}$ because $x^{l^n}\in U_i$ and $F=Z(W)$ for a Wedderburn component $W$ of $\Q_l(\zeta_i)\otimes_{\Q_l}\QQ U_i$.)
By Corollary \ref{Cor1}, together with the precondition that 
$SK_1(\Q_l(\zeta_i)\otimes_{\Q_l}\QQ U_i)=1$ and in particular $SK_1(W)=1$, this also implies
\begin{Prop}
With the above notations, assume $W^{x^{j}}\neq W$ for all $1\leq j\leq l^n-1$ and $SK_1(\Q_l(\zeta_i)\otimes_{\Q_l}\QQ U_i)=1$ for all $i$. Then
$$SK_1(\tilde{W})=SK_1(W)=1.$$ 
\end{Prop}
\hfill$\Box$

\subsubsection{The case $d=0$}

Next, we consider $d=0$, i.e.~$W^x=W$.  Thus,
$$\tilde{W}=\bigoplus_{j=0}^{l^n-1} Wx^j,\ \ V=\tilde{W}.$$ This time, $Z(\tilde{W})=F^{\langle x\rangle}=E$ with $[F:E]=l^n$ and $G(F/E)=\langle \sigma\rangle$, where $\sigma$ is induced by the conjugation by $x$. Note that 
$$\langle \sigma\rangle \cong\langle\overline{x}\rangle\cong\langle\tau\rangle\subseteq G(\Q_l(\zeta_i)/\Q_l)$$
with $\tau$ as above induced by the action of $x$ on $\Q_l(\zeta_i)$.

First, we give a brief outline of the proof of $SK_1(\tilde{W})=SK_1(V)=1$.
\setcounter{St}{0}
\begin{St}
$F\otimes_E V\cong W_{l^n\times l^n}\subseteq V_{l^n\times l^n}$.
\end{St}
\begin{St}
There exists a $w\in W_{l^n\times l^n}$ such that the conjugation by $w^{-1}x$ is the automorphism $C_{w^{-1}x}=\sigma\otimes 1$  on $W_{l^n\times l^n}$. It is of order  $l^n$ and $(w^{-1}x)^{l^n}\in Z(W_{l^n\times l^n})=F$.
\end{St}
\begin{St}
$(w^{-1}x)^{l^n}=1$ and therefore $A=(F/E, \sigma, (w^{-1}x)^{l^n})\subseteq V_{l^n\times l^n}$
is a central simple split $E$-algebra.
\end{St}
\begin{St}
$Z_{V_{l^n\times l^n}}(A)=(W_{l^n\times l^n})^{\langle w^{-1}x\rangle}$.
\end{St}
\begin{St}
$V_{l^n\times l^n}\cong A\otimes_E (W_{l^n\times l^n})^{\langle w^{-1}x\rangle}$.
\end{St}
\begin{St}
$SK_1(V)=SK_1((W_{l^n\times l^n})^{\langle w^{-1}x\rangle})=1$.
\end{St}

We now start with the sketched computations.
We can read $V$ as free left $W$-module of rank $l^n$ with basis $1,x,...,x^{l^n-1}$, i.e. $V=\bigoplus_{j=0}^{l^n-1}W x^j$.  This allows us to formulate

\begin{lemma}\label{isom}
With the above notations, we have an isomorphism
$$F\otimes_E V\stackrel{\cong}{\longrightarrow} W_{l^n\times l^n}=\Hom_W(V,V),\ \ f\otimes v\mapsto l_f\circ r_v,$$
where $l_f$ resp.~$r_v$ denotes the left resp.~right multiplication with $f\in F$ resp.~$v\in V$.
\end{lemma}
\begin{remark}
In particular, $f\otimes 1\in F\otimes_E V$ maps to the diagonal matrix $f\cdot \textbf{1}$ with $\textbf{1}$ the unit matrix.
\end{remark}

\textbf{Proof:}
$F\otimes_E V$ and $W_{l^n\times l^n}$ are isomorphic by \cite[Cor 7.14]{Re}, we only have to substitute $K$ by $E$, $A$ by $V$ and $B$ by $F$. Then, we get $r=l^n$ and the centralizer  $B'=Z_V(F)=W$ implies $F\otimes_E V=F^{op}\otimes_E V\cong W_{l^n\times l^n}$.

We will call the stated isomorphism $\varphi$ for the moment. We read the actions of the $W$-endomorphisms of $V$ by the right to ensure that $\varphi$ is compatible with multiplication. For this, take $f,\ f'\in F$ and $v,\ v',\ a\in V$.  Then the commutativity of $F$ yields
\begin{align*}
(a)(\varphi(f\otimes v)\circ\varphi(f'\otimes v'))&=(a)((l_f\circ r_v)\circ (l_{f'}\circ r_{v'}))\\
&=f'favv'=ff'avv'\\
&=(a)(l_{ff'}\circ r_{vv'})=(a)\varphi(ff'\otimes vv')\\
&=(a)\varphi((f\otimes v)(f'\otimes v')).
\end{align*}

It now easily follows that $\varphi$ is a homomorphism of $E$-algebras.

$F\otimes_E V$ is simple because $V$ is a central simple $E$-algebra. Thus, $\varphi\neq 0$ implies that $\varphi$ is injective. By dimension comparison, it is surjective as well. \hfill$\Box$

Next, we construct the automorphism $C_{w^{-1}x}$ on $W_{l^n\times l^n}$.
On the one hand,  conjugation by $x$ is an automorphism $c_x$ on $W$ and can therefore be extended to $W_{l^n\times l^n}$ by letting it act on the matrix entries. Furthermore, we can read $x$ as the diagonal matrix $M_x=x\cdot \textbf{1}$ in $V_{l^n\times l^n}$. Then, the extension of $c_x$ on $W_{l^n\times l^n}$ is the conjugation by this matrix $M_x$. This automorphism on $W_{l^n\times l^n}$ will be called $C_x$ in the sequel and we remark that $C_x$ acts on $F=F\cdot \textbf{1}=F\otimes_E E$ as $\sigma$, with $\langle \sigma\rangle=G(F/E)$ as above.

On the other hand, $\sigma\otimes 1:\ F\otimes_E V\rightarrow F\otimes_E V$ is another automorphism on $W_{l^n\times l^n}$. As the restriction of $\sigma\otimes 1$ to $F\otimes_E E$ is by construction the old isomorphism $\sigma$, the actions of $C_x$ and $\sigma\otimes 1$  coincide on $F=F\otimes_E E$. 

Therefore, $C_x(\sigma\otimes 1)^{-1}$ is a central automorphism on $W_{l^n\times l^n}$, i.e.~it acts trivially on the centre $Z(W_{l^n\times l^n})=F\cdot \textbf{1}=F$. The theorem of Skolem-Noether now implies that $C_x(\sigma\otimes 1)^{-1}$ is the conjugation $C_w$ by some $w\in W_{l^n\times l^n}$, i.e.
$$\sigma \otimes 1=C_w^{-1}C_x=C_{w^{-1}x}.$$
As $(\sigma\otimes 1)^{l^n}=\id$, we furthermore conclude $(C_{w^{-1}x})^{l^n}=\id$. This means that the conjugation by $$(w^{-1}x)^{l^n}=(w^{-1})^{1+x^{-1}+...+x^{-l^n+1}}x^{l^n}$$
is trivial on $W_{l^n\times l^n}$ and, as $x^{l^n}\in W$ (more precisely $x^{l^n}\in \Q(\zeta_i)\otimes_{\Q_l}\QQ U_i$ has a component in $W$ but we suppress  this here for the sake of brevity), we conclude
$$(w^{-1}x)^{l^n}\in Z(W_{l^n\times l^n})=F\cdot \textbf{1}=F.$$

Finally, we choose
$$A=(F\otimes_E E/E\otimes_E E,\sigma\otimes 1=C_{w^{-1}x},(w^{-1}x)^{l^n})=(F/E,\sigma, (w^{-1}x)^{l^n}).$$
By  construction, we have $w\in W_{l^n\times l^n}$ and $x\in V$. Therefore, $$w^{-1}x\in \bigoplus_{j=0}^{l^n-1}W_{l^n\times l^n}x^j=V_{l^n\times l^n}$$
 and hence $A\subseteq V_{l^n\times l^n}$.

\begin{lemma}\label{ll1}
Let $A=(F/E,\sigma, (w^{-1}x)^{l^n})$ be as above. Then $A$ splits, i.e.~$A\sim E$ in $\Br(E)$.
\end{lemma}
\textbf{Proof:}
The cyclic algebra $A$ splits if $(w^{-1}x)^{l^n}$ is a norm element in $E$, i.e.~if there exists an element $f\in F$ with $\NN_{F/E}(f)=(w^{-1}x)^{l^n}$, where $\NN_{F/E}=\NN_{\langle \sigma \rangle}$ is the Galois norm of the field extension $F/E$.
To show this, we compute $(w^{-1}x)^{l^n}$ explicitely. 

First, let $k_x$ denote the conjugation by $x$ on $V$, s.t. we can study the automorphism 
$$\sigma\otimes k_x:F\otimes_E V\rightarrow F\otimes_E V.$$
By Lemma \ref{isom}, we know
$$F\otimes_E V\cong W_{l^n\times l^n}=\Hom_W(V,V),\ \ f\otimes v\mapsto l_f\circ r_v.$$
 Here, we choose a basis $1, x, ..., x^{l^n-1}$ of the free left $W$-vector space $V$. Then, we write $v=\sum_{i=0}^{l^n-1}w_ix^i$ and achieve that $l_f\circ r_v$ is represented by the matrix
$$\left( \begin{array}{cccc}
fw_0 & fw_1 & \dots & fw_{l^n-1}\\
fw_{l^n-1}^{x{^{-1}}}x^{l^n} & fw_0^{x^{-1}} & \dots  & . \\
\vdots & \vdots & \vdots & \vdots \\
fw_1^{x^{-(l^n-1)}}x^{l^n} &. &\dots &fw_0^{x^{-(l^n-1)}}
\end{array}\right).$$
Here, we recall that we write the matrices from the right. Next,
$$(\sigma\otimes k_x)(f\otimes v)=\sigma(f)\otimes v^x \leftrightarrow
\left( \begin{array}{cccc}
\sigma(f)w_0^x & \sigma(f)w_1^x & \dots &\sigma(f)w_{l^n-1}^x \\
\sigma(f)w_{l^n-1}x^{l^n} &\sigma(f)w_0 & \dots  & . \\
\vdots & \vdots & \vdots & \vdots \\
\sigma(f)w_1^{x^{-(l^n-2)}}x^{l^n} &.&\dots & .
\end{array}\right).$$

A comparison of the two matrices shows that $\sigma\otimes k_x$ is the conjugation by $x$ on $W_{l^n\times l^n}$, i.e. 
$$\sigma\otimes k_x =C_x.$$

Next, we obtain
$$C_x\circ(\sigma\otimes 1)^{-1} =(\sigma\otimes k_x)\circ(\sigma\otimes1)^{-1}=1\otimes k_x=C_{1\otimes x}:\ F\otimes_E V\rightarrow F\otimes_E V.$$

Now, we see that 
$$1\otimes x=w\in W_{l^n\times l^n}$$
is the element s.t. $C_{w^{-1}x}=\sigma\otimes 1$.

Read as matrices in $V_{l^n\times l^n}$, we can write
\begin{align*}
w^{-1}x=(x^{-1}w)^{-1}&=\left(\left(
\begin{array}{ccccc}
x^{-1} & \\
&x^{-1} \\
&&\ddots\\
&&&\ddots\\
&&&&x^{-1}
\end{array}\right)
\cdot\left(
\begin{array}{ccccc}
0&1&0&\dots & 0\\
0&0&1&\dots & 0\\
\vdots &\vdots& \vdots &\vdots&\vdots\\
0&0&0&\dots &1\\
x^{l^n}&0&0&\dots&0
\end{array}\right)\right)^{-1}
\\
&=\left(
\begin{array}{ccccc}
0&x^{-1}&0&\dots &0\\
0 &0 &x^{-1}&\dots&0\\
\vdots & \vdots &\vdots &\vdots &\vdots\\
0 &0  & 0&\dots&x^{-1}  \\
x^{l^n-1}&0&0&\dots &0
\end{array}\right)^{-1}
\end{align*}

We finally conclude that 
$$(w^{-1}x)^{l^n}=1$$
which is certainly a norm element.
\hfill$\Box$

\begin{lemma}
Set $A=(F/E,\sigma, (w^{-1}x)^{l^n})$ and $V$ as above. Then
$$Z_{V_{l^n\times l^n}}(A)= (W_{l^n\times l^n})^{\langle w^{-1}x\rangle}.$$
\end{lemma}
\textbf{Proof:}
First, we have to read $A$ in $V_{l^n\times l^n}$.
For this, we observe that $E$ and $F$ are to be represented by the diagonal matrices $E=E\cdot\textbf{1}$ and $F=F\cdot \textbf{1}$. Then, choose a matrix $(v_{ij})_{i,j}\in Z_{V_{l^n\times l^n}}(A)$ and $f\cdot\textbf{1}\in F\cdot \textbf{1}$. We get
$$f^{-1}\textbf{1}(v_{ij})f\textbf{1}=(f^{-1}v_{ij}f)\stackrel{!}{=}(v_{ij}),$$
i.e.~$f^{-1}v_{ij}f\stackrel{!}{=}v_{ij}$ for all $i,j=0,...,l^n-1$. As this equation has to be fulfilled for all $f\in F$, but $v_{ij}=w_0+...+w_{l^n-1}x^{l^n-1}\in V$, we conclude that $v_{ij}=w_0\in W$. Therefore, $Z_{V_{l^n\times l^n}}(A)\subseteq W_{l^n\times l^n}$. 

Next, we conjugate by $w^{-1}x$ and see
$Z_{V_{l^n\times l^n}}(A)\subseteq (W_{l^n\times l^n})^{\langle w^{-1}x\rangle}$. 

 It remains to show the other inclusion $(W_{l^n\times l^n})^{\langle w^{-1}x\rangle}\subseteq Z_{V_{l^n\times l^n}}(A)$. But this is obvious, because $A=\bigoplus_{j=0}^{l^n-1} F\cdot \textbf{1}(w^{-1}x)^{j}$ and $v\in (W_{l^n\times l^n})^{\langle w^{-1} x\rangle}$ commutes with $w^{-1}x$ as well as with $a\in F\cdot \textbf{1}=Z(W_{l^n\times l^n})$. Thus, $v$ commutes with $a_0+...+a_{l^n-1}(w^{-1}x)^{l^n-1}\in \bigoplus_{j=0}^{l^n-1}F(w^{-1}x)^{j}=A$, too.
\hfill$\Box$

\begin{corollary} \label{lem123}
$(W_{l^n\times l^n})^{\langle w^{-1}x\rangle}$ is a central simple $Z(A)=E$-algebra.
\end{corollary}
\textbf{Proof:}
This is true by the centralizer theorem. \hfill$\Box$

\begin{lemma}\label{ll2}
With the above notations, we have
$$V\cong  (W_{l^n\times l^n})^{\langle w^{-1}x\rangle}.$$
Moreover,
$$ F\otimes_E (W_{l^n\times l^n})^{\langle w^{-1}x\rangle} \stackrel{\cong}{\longrightarrow} W_{l^n\times l^n},\ \ f\otimes w\mapsto fw,$$
and 
$$ A\otimes_E Z_{V_{l^n\times l^n}}(A) \stackrel{\cong}{\longrightarrow} V_{l^n\times l^n},\ \ a\otimes v\mapsto av,$$
are isomorphisms.
\end{lemma}
\textbf{Proof:}
First, $V\sim  (W_{l^n\times l^n})^{\langle w^{-1}x\rangle}$ in $\Br(E)$ by the centralizer theorem which states that
$$Z_{V_{l^n\times l^n}}(A)\sim A^{op}\otimes_E V_{l^n\times l^n}$$
in $\Br(E)$. By Lemma \ref{ll1}, we know that $A\sim E$ in $\Br(E)$. Because $A^{op}$ is the inverse of $A$ in $\Br(E)$, we conclude $A^{op}\sim E$ in $\Br(E)$, too. Thus,
$$(W_{l^n\times l^n})^{\langle w^{-1}x\rangle}=Z_{V_{l^n\times l^n}}(A)\sim A^{op}\otimes_E V_{l^n\times l^n}\sim E\otimes_E V_{l^n\times l^n}\sim V.$$

Next, we compute the respective degrees over $E$:
\begin{align*}
[(W_{l^n\times l^n})^{\langle w^{-1}x\rangle}:E]=[W_{l^n\times l^n}:E]/l^n=l^n[W:E]=l^n[W:F][F:E]=l^{2n}[W:F]
\end{align*}
and
\begin{align*}
[V:E]=[V:W][W:F][F:E]=l^{2n}[W:F].
\end{align*}
Thus, $V$ and $(W_{l^n\times l^n})^{\langle w^{-1}x\rangle}$ are as Brauer equivalent algebras of the same degree isomorphic.

We turn to the second isomorphism. As $Z_{V_{l^n\times l^n}}(A)=(W_{l^n\times l^n})^{\langle w^{-1} x\rangle}$ is a central simple $E$-algebra, $F\otimes_E (W_{l^n\times l^n})^{\langle w^{-1}x\rangle}$ is a central simple $F$-algebra. We compute the respective degrees over $F$:
$$[V_{l^n\times l^n}:E]=[V_{l^n\times l^n}:W_{l^n\times l^n}][W_{l^n\times l^n}:F][F:E]=l^{2n}[W_{l^n\times l^n}:F]$$
and, by the centralizer theorem,
\begin{align*}
[V_{l^n\times l^n}:E]&=[A:E][Z_{V_{l^n\times l^n}}(A):E]\\
&=[F:E]^2[(W_{l^n\times l^n})^{\langle w^{-1}x\rangle} :E]=l^{2n}[(W_{l^n\times l^n})^{\langle w^{-1}x \rangle}:E]
\end{align*}
imply 
\begin{align*}
[W_{l^n\times l^n}:F]&=[(W_{l^n\times l^n})^{\langle w^{-1}x\rangle}:E]=[(W_{l^n\times l^n})^{\langle w^{-1}x\rangle}\otimes_E F:E\otimes_E F]\\
&=[(W_{l^n\times l^n})^{\langle w^{-1} x\rangle}\otimes_E F:F].
\end{align*}
 Next,  $F\otimes_E (W_{l^n\times l^n})^{\langle w^{-1}x\rangle}\rightarrow W_{l^n\times l^n}$, $f\otimes w\mapsto fw$, is injective because otherwise the kernel would form a non-trivial two-sided ideal. But  $F\otimes_E (W_{l^n\times l^n})^{\langle w^{-1}x\rangle}$ is a central simple $F$-algebra. Thus, the only non-trivial two-sided ideal is $F\otimes_E (W_{l^n\times l^n})^{\langle w^{-1}x\rangle}$ itself, which is impossible because $E \otimes_E (W_{l^n\times l^n})^{\langle w^{-1}x\rangle}\subseteq F\otimes_E (W_{l^n\times l^n})^{\langle w^{-1}x\rangle}$ maps to $(W_{l^n\times l^n})^{\langle w^{-1}x \rangle}\subseteq W_{l^n\times l^n}$ and thus the kernel can not be $F\otimes_E (W_{l^n\times l^n})^{\langle w^{-1}x\rangle}$.
This implies $F\otimes_E (W_{l^n\times l^n})^{\langle w^{-1}x\rangle}\subseteq W_{l^n\times l^n}$. As both sides are of the same degree over $F$, we conclude $F\otimes_E (W_{l^n\times l^n})^{\langle w^{-1}x\rangle}=W_{l^n\times l^n}$.

Finally, we show $A\otimes_E Z_{V_{l^n\times l^n}}(A)\cong V_{l^n\times l^n}$. As $A$ and $Z_{V_{l^n\times l^n}}(A)$ are central simple $E$-algebras, $A\otimes_E Z_{V_{l^n\times l^n}}(A)$ is a central simple $E$-algebra, too. We again show that the respective degrees over $E$ coincide:
$$[V_{l^n\times l^n}:E]=[A:E][Z_{V_{l^n\times l^n}}(A):E]=[A\otimes_E Z_{V_{l^n\times l^n}}(A):E].$$
Next, the homomorphism $A\otimes_E Z_{V_{l^n\times l^n}}(A)\rightarrow V_{l^n\times l^n}$, $a\otimes v\mapsto av$, again is injective. Dimension comparison implies $A\otimes_E Z_{V_{l^n\times l^n}}(A)=V_{l^n\times l^n}$.
\hfill$\Box$

\begin{Prop}
With the above notations,  assume that $W^{x}=W$ and moreover
\\${SK_1(\Q_l(\zeta_i)\otimes_{\Q_l}\QQ U_i)=1}$ for all $i$. Then
$$SK_1(\tilde{W})=SK_1((W_{l^n\times l^n})^{\langle w^{-1}x \rangle})=1.$$
\end{Prop}
\textbf{Proof:}
We are still in the case $V=\tilde{W}$. With   $V_{l^n\times l^n}= A\otimes_E Z_{V_{l^n\times l^n}}(A)$, it therefore suffices to compute 
\begin{align*}
SK_1(V)&=SK_1(V_{l^n\times l^n})= SK_1(A\otimes_E (W_{l^n\times l^n})^{\langle w^{-1}x\rangle})\\
&\stackrel{1}{=} SK_1(A)\times SK_1((W_{l^n\times l^n})^{\langle w^{-1}x\rangle})\\
&\stackrel{2}{=}1\times SK_1((W_{l^n\times l^n})^{\langle w^{-1}x\rangle}).
\end{align*}
in the sequel. For  $\stackrel{2}{=}$, we use that $A$ splits and therefore 
 $SK_1(A)=1$; moreover, $A$ and $(W_{l^n\times l^n})^{\langle w^{-1}x\rangle}$ have coprime Schur indices which implies $\stackrel{1}{=}$ by \cite[Lem 5, p.~160]{Drax}.

Now, we choose a $v\in V_{l^n\times l^n}$ with $\nr_{V_{l^n\times l^n}/E}(v)=1$. It represents an element in 
$SK_1(V_{l^n\times l^n})$.
By the above, the class of $v$  can be read as 
$$[v]=(1, [\tilde{v}])=[1\otimes \tilde{v}]=[\tilde{v}]$$
  with $\tilde{v}\in (W_{l^n\times l^n})^{\langle w^{-1}x\rangle}\subseteq V_{l^n\times l^n}$ and 
$$\nr_{(W_{l^n\times l^n})^{\langle w^{-1}x\rangle}/E}(\tilde{v})=1.$$

Therefore,  $v$ and $\tilde{v}$  only differ by a factor in $[(V_{l^n\times l^n})^\times,(V_{l^n\times l^n})^\times]$.
It hence suffices to show that $\tilde{v}\in [(V_{l^n\times l^n})^\times,(V_{l^n\times l^n})^\times]$ for $SK_1(V_{l^n\times l^n})=1$.

For the computation of $\nr_{(W_{l^n\times l^n})^{\langle w^{-1}x\rangle}/E}$, let $M$ be a splitting field of $(W_{l^n\times l^n})^{\langle w^{-1}x\rangle}$ with $M\supseteq F$. Thus, as 
$$(W_{l^n\times l^n})^{\langle w^{-1}x\rangle}\subseteq F\otimes_E  (W_{l^n\times l^n})^{\langle w^{-1}x\rangle} =W_{l^n\times l^n},$$
 we get
$$M_{m\times m}= M\otimes_E (W_{l^n\times l^n})^{\langle w^{-1} x\rangle}=M\otimes_F F\otimes_E (W_{l^n\times l^n})^{\langle w^{-1}x \rangle}=M\otimes_F W_{l^n\times l^n}$$
for a certain $m\in \N$, i.e.~$M$ is also a splitting field of $W_{l^n\times l^n}$. This implies
$$1=\nr_{(W_{l^n\times l^n})^{\langle w^{-1}x\rangle}/E}(\tilde{v})\stackrel{1}{=}\nr_{W_{l^n\times l^n}/F}(\tilde{v}),$$
where $\stackrel{1}{=}$ holds due to the isomorphism $F\otimes_E (W_{l^n\times l^n})^{\langle w^{-1}x \rangle}=W_{l^n\times l^n}$, $1\otimes \tilde{v}\mapsto 1\cdot \tilde{v}$ and the common splitting field $M\supseteq F\supseteq E$ of $(W_{l^n\times l^n})^{\langle w^{-1}x \rangle}$ and $W_{l^n\times l^n}$.
But, by assumption, $SK_1(W_{l^n\times l^n})=SK_1(W)=1$ and hence
$$\tilde{v}\in [(W_{l^n\times l^n})^\times,(W_{l^n\times l^n})^\times]\subseteq [(V_{l^n\times l^n})^\times,(V_{l^n\times l^n})^\times].$$
This concludes the proof.
\hfill$\Box$

\subsubsection{The intermediate case $0<d<n$}

Finally, the triviality of $SK_1(\tilde{W})$ in the intermediate cases for $0<j<n$ is a consequence of the extreme cases: We fix a $0<d<n$ and set $y:=x^{l^d}$ and $m:=n-d$. Thus, 
\begin{align*}
V&=W\oplus Wx^{l^d}\oplus...\oplus Wx^{l^d(l^{n-d}-1)}\\
&=W\oplus Wy\oplus...\oplus Wy^{(l^{n-d}-1)}=\bigoplus_{j=0}^{l^m-1}Wy^j.
\end{align*}
As $\tilde{W}=V_{l^d\times l^d}$, it suffices to compute $SK_1(V)$. But $V$ is now of the same form as $\tilde{W}$ in the case $d=0$, with $x$ replaced by $y$ and $n$ replaced by $m$. Thus, we only have to check that the above arguments apply to this $V$ in the same manner. As it can be seen easily that we can copy the above literally, we leave this to the reader.

Hence, we have seen that $SK_1(\tilde{W})=1$ for every Wedderburn component of $\QQ G$. This concludes the proof of Theorem \ref{sk1}. \hfill$\Box$

\subsection{$\Q_l$-$q$-elementary groups $G$}\label{l-q-section}

Next, we consider the case of $\Q_l$-$q$-elementary groups $G$ with $q\neq l$. Here, our result on the triviality of the reduced Whitehead group is stronger than in the case $q=l$ because it holds without assumptions:

\begin{theorem}\label{sk2}
Let $G$ be a $\Q_l$-$q$-elementary group with $q\neq l$ prime. Then 
 $$SK_1(\QQ G)=1.$$
\end{theorem}

The proof of this theorem closely follows the proof of Theorem \ref{sk1}. Thus, we only give a short outline how to adapt the ideas used for the case $q=l$ to our new situation.

To do so, we first recall that the $\Q_l$-$q$-elementary group $G$ is a direct product $G=H\times \Gamma$ with $H$ a finite $\Q_l$-$q$-elementary group. More precisely, $H=\langle s \rangle \rtimes H_q$ with $\langle s \rangle$ a cyclic group of order prime to $q$ and a $q$-group $H_q$ whose action on $\langle s \rangle$ induces a homomorphism $H_q \rightarrow G(\Q_l(\zeta)/\Q_l)$ for $\zeta$ a primitive root of unity of order $|\langle s\rangle|$.

We now take $\langle s_l\rangle$ the $l$-Sylow subgroup of $\langle s\rangle$, thus $\langle s\rangle = \langle s_l\rangle\times\langle s'\rangle$, and obtain $G=\langle s_l\rangle \rtimes U$ with $U=(\langle s'\rangle \rtimes H_q)\times \Gamma$ and $\langle s'\rangle \rtimes H_q$ an $l$-prime group. Still, $U$ acts on $\langle s_l\rangle$ via Galois automorphisms. 

As in the $\Q_l$-$l$-elementary case, we fix a finite set $\{\beta_i\}$ of representatives of the $G(\Q_l^c/\Q_l)$-orbits of the irreducible $\Q_l^c$-characters of $\langle s_l \rangle$. Let also $\zeta_i$ denote a fixed primitive root of unity with $\beta_i(s)=\zeta_i$.  This time, $\Q_l(\zeta_i)$ is not unramified because $\langle s_l\rangle$ is an $l$-group, but $\Q_l(\zeta_i)/\Q_l$ remains cyclic because $l$ is odd. 
Let again $U_i:=\{u\in U: \beta_i^{u}=\beta_i\}$ denote the stabilizer group of $\beta_i$.  Thus, $U_i\lhd U$ and $A_i:=U/U_i\leq G(\Q_l(\beta_i)/\Q_l)=G(\Q_l(\zeta_i)/\Q_l)$ is cyclic because $G(\Q_l(\zeta_i)/\Q_l)$ is. We fix a representative $x\in U$ with $\langle \overline{x}\rangle =U/U_i=A_i$.  Then, $\overline{x}$ maps to some $\tau$ under the injection $U/U_i\rightarrowtail G(\Q_l(\zeta_i)/\Q_l)$; and  $|x|=|U/U_i|=:l^n$ with again $x$, $\tau$ and $n$ depending on $i$.
Finally, we set $G_i:=\langle s_l \rangle \rtimes U_i$.

Now, we may compute the isomorphism
$$\QQ G=\bigoplus_i e_i\QQ G\cong \bigoplus_i (\Q_l(\zeta_i)\otimes_{\Q_l}\QQ U_i)\star \langle x \rangle.$$

In this case that $q\neq l$, the group $U_i$ is not pro-$l$ and therefore the structure of $\QQ U_i$ differs from the structure of the analogous object in the pro-$l$ case as stated in Section \ref{Recollections}. Yet, we may collect all relevant information on $\QQ U_i$ easily. Recall that $U$ is the direct product $U=(\langle s'\rangle\rtimes H_q)\times \Gamma$ and therefore $U_i$ also is the direct product $U_i=H'\times \Gamma$ with $H'$ a subgroup of $\langle s'\rangle \rtimes H_q$. Now, \cite[74.11, p. 740]{CR2} implies that $\QQ U_i$ is the direct sum of matrix rings over the fields $F'=\Q_l(\eta')\otimes_{\Q_l}\QQ\Gamma$ for certain characters $\eta'$ of $H'$. Because $F'$ and $\Q_l(\zeta_i)$ are linearly disjoint over $\Q_l$, the algebra $\Q_l(\zeta_i)\otimes_{\Q_l}\QQ U_i$ also is a direct sum of matrix rings over fields $F=\Q_l(\zeta_i)\otimes_{\Q_l}F'$. This proves 

\begin{Prop}
With the above notations, we have
$$SK_1(\Q_l(\zeta_i)\otimes_{\Q_l}\QQ(U_i))=1.$$
\end{Prop}
This proposition allows us to formulate Theorem \ref{sk2} in the stronger form without any assumptions on $SK_1(\Q_l(\zeta_i)\otimes_{\Q_l}\QQ(U_i))$.

It now remains to examine the semisimple algebra $(\Q_l(\zeta_i)\otimes_{\Q_l}\QQ U_i)\star \langle x\rangle$. We may transfer the computation for the case $q=l$ directly to our new situation.

This concludes the proof of Theorem \ref{sk2}. \hfill$\Box$

\bibliographystyle{abbrv}
\bibliography{literature_lau}

\end{document}